\documentclass[12pt]{amsart}  
\usepackage[T1]{fontenc}
\usepackage{amsxtra}

\usepackage{color}
\usepackage{amscd,amsfonts,amsmath,amssymb,
amsthm,amstext,latexsym}
\usepackage{enumerate,pifont,graphicx}
\usepackage[english]{babel}
\usepackage[all,cmtip]{xy} 
\newcommand{\Z}{{\mathbb Z}}

\newcommand{\R}{{\mathbb R}}
\newcommand{\C}{{\mathbb C}}

\def\boC{{\mathcal C}}

\def\cqfd{\hfill$\Box$}
\def\Res{{\,\rm Res}}

\def\Re{{\rm Re}}
\def\Im{{\rm Im}}

\def\whOmega{{\widehat{\Omega}}}
\def\whSigma{{\widehat{\Sigma}}}

\def\whgamma{{\widehat{\gamma}}}
\def\ii{{\rm i}}

\newtheorem{theorem}{Theorem}
\newtheorem{lemma}{Lemma}
\newtheorem{proposition}{Proposition}
\newtheorem{remark}{Remark}

\newtheorem{definition}{Definition}
\newtheorem{hypothesis}{Hypothesis}

\title{Hollow vortices and minimal surfaces}

\author{Martin Traizet}
\thanks{Partially supported by ANR-11-ISO1-0002 grant.}

\begin{document}
\maketitle

\bigskip

\bigskip
\section{Introduction}
\label{section1}
We consider the following overdetermined problem in the plane, known as the
 {\em hollow vortex problem}:
\begin{equation}
\label{HVP}
\left\{\begin{array}{ll}
\Delta u =0 &\mbox{ in $\Omega$}\\
u \mbox{ constant } & \mbox{ on each component of $\partial\Omega$}\\
||\nabla u||=1& \mbox { on $\partial\Omega$}\\
\end{array}\right.
\end{equation}
Here $\Omega\subset\R^2$ is an unbounded domain with smooth, non-empty boundary, and
$u:\Omega\to\R$ is a smooth function. Observe that we require $u$ to be constant
on each boundary component, but we do not ask that the constant is the same for
all boundary components.
The fact that $||\nabla u||=1$ on the boundary is equivalent to the Neumann condition
$\frac{\partial u}{\partial\nu}=\pm 1$, where $\nu$ is the (interior) unit normal to the boundary.
Again, the sign of $\frac{\partial u}{\partial\nu}$ may depend on the boundary component.

\begin{figure}
\begin{center}
\includegraphics[height=40mm]{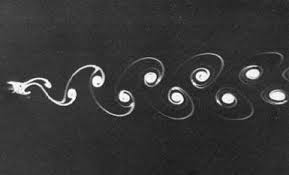}
\end{center}
\caption{A Von Karman vortex street, from Van Dyke Album of Fluid Motion
\cite{vandyke}}
\label{fig1}
\end{figure}

\medskip
Problem \eqref{HVP} is overdetermined because both Neumann and Dirichlet
boundary conditions are prescribed, which is not possible for general
domains $\Omega$.
A domain $\Omega$ admitting a function $u$ solving Problem \eqref{HVP} will
be called a {\em domain with hollow vortices}. The hollow vortices refer
to the components of $\R^2\setminus\Omega$. The name comes from
the following physical interpretation of Problem \eqref{HVP}:
the stationary flow of an inviscid, incompressible fluid in a domain $\Omega$ is described by Euler 
equations:
$$\mbox{div } \vec{v}=0,\qquad
(\vec{v}\cdot\nabla)\vec{v}=-\frac{1}{\rho}\nabla p$$
where $\vec{v}$ denotes  the velocity vector, $p$ the pressure and $\rho$ the mass density
of the fluid.
In the 2-dimensional case, we can write $\vec{v}=(\frac{\partial u}{\partial y},
-\frac{\partial u}{\partial x})$ for some function $u$ called the stream function.
If we assume moreover that the flow is irrotational, then $\Delta u=0$.
The condition that $u$ is constant on a boundary component $\gamma$ of $\Omega$ means that
$\gamma$ is a stream line. The condition that $||\nabla u||$
is constant on $\gamma$ means that the norm of the velocity is constant,
which is equivalent to the fact that the pressure is constant by
Bernoulli law. We can think of $\gamma$ as bounding a spinning bubble of air
with constant pressure inside, or ``hollow vortex''.
Note that we require the pressure to be the same on each boundary component,
so we should call \eqref{HVP} the {\em isobaric} hollow vortex problem.
Hollow vortices have been proposed as a model for some periodic configurations of
vortices observed in the turbulent flow past an obstacle known as Von Karman vortex streets
\cite{crowdy}, see Figure \ref{fig1}.

\medskip
In the particular case where $u=0$ on $\partial\Omega$ and $u>0$ in
$\Omega$, solutions to the hollow vortex problem have been studied in
\cite{hhp}, \cite{khavinson} and completely classified by the author
in \cite{GAFA}, by establishing a correspondence with a certain type
of minimal surfaces and using classification results in minimal surface theory.
It turns out that the correspondence extends to the general case of Problem
\eqref{HVP}, only to a wider class of minimal surfaces.
Our goal in this paper is to describe this correspondence and give examples.
\medskip

The corresponding overdetermined problem for minimal surfaces is the following:
\begin{equation}
\label{mse}
\left\{\begin{array}{ll}
(1+v_y^2)v_{xx}+(1+v_x^2)v_{yy}-2v_x v_y v_{xy}=0 &\mbox{ in $\whOmega$}\\
v\mbox{ constant }&\mbox{ on each component of $\partial\whOmega$}\\
||\nabla v||\to\infty &\mbox{ on $\partial\whOmega$}\\
\end{array}\right.
\end{equation}
Here $\whOmega\subset\R^2$ is an unbounded domain with non-empty boundary
and $v:\whOmega\to\R$ is a smooth function. The subscripts denote partial
derivatives. The first equation is the minimal surface equation:
it says that the graph of $v$, denoted $M$, is a minimal surface.
The limit in the last condition means the following:
for any $z_0\in\partial\whOmega$, $\lim_{z\to z_0} \frac{\partial v}{\partial \nu}=\pm\infty$, the
limit being uniform on compact sets of $\partial\whOmega$.
Geometrically speaking, this means that the Gauss map (i.e. the unit normal vector) of $M$ is horizontal on the boundary.
Such minimal surfaces can be smoothly extended by reflection in the horizontal plane
containing each boundary component. For this reason, we call them {\em minimal
graphs bounded by horizontal symmetry curves}.
\medskip

For simplicity, we assume that all domains $\Omega$ and $\whOmega$ considered in this paper
satisfy the following
finiteness hypothesis:
\begin{hypothesis}
\label{finiteness}
Either:
\begin{itemize}
\item $\Omega$ has a finite number of boundary components,
\item or $\Omega$ is invariant by a translation $T$ and the quotient
$\Omega/T$ has a finite number of boundary components
(simply-periodic case),
\item or $\Omega$ is invariant by two independent translations
(doubly-periodic case).
\end{itemize}
\end{hypothesis}
Under this hypothesis, the main result of this paper is
\begin{theorem}
\label{theorem1}
There is a 1:1 correspondence between:
\begin{itemize}
\item solutions $(\Omega,u)$ of Problem \eqref{HVP} such that
$||\nabla u||<1$ in $\Omega$,
\item solutions $(\whOmega,v)$ of Problem \eqref{mse}.
\end{itemize}
\end{theorem}
We describe how the correspondence works in Section \ref{section2}.
An interesting feature of the correspondence is that the domain with hollow vortices
$\Omega$ and
the corresponding minimal graph $M$ are conformally related. Also,
each component of $\partial\whOmega$ is a translation of the corresponding
component of $\partial\Omega$.

\medskip
If $\Omega$ is a domain with hollow vortices, then in general,
extending the corresponding minimal surface $M$ by reflection will not give an embedded
minimal surface.
However, there are two particular cases where $M$ can be extended to a complete embedded minimal surface:
\begin{itemize}
\item Case I: $u>0$ in $\Omega$ and $u=0$ on $\partial \Omega$
(this is the case
already considered in \cite{GAFA}). In this case, the corresponding minimal surface
$M$ lies in the half-space $x_3>0$ and has boundary in the horizontal plane
$x_3=0$. It can be extended by reflection to a complete, embedded minimal surface.
\item Case II: $0<u<c$ in $\Omega$ and $u=0$ or $u=c$ on each boundary
component of $\partial\Omega$. In this case, $M$
lies in the slab $0<x_3<c$ and has boundary in the horizontal planes at height
$0$ and $c$. It can be extended by iterated reflections in horizontal planes into
a complete, embedded, periodic minimal surface with vertical period $2c$.
\end{itemize}
We know a lot of such minimal surfaces, including some triply-periodic minimal
surfaces discovered in the 19th century by Schwarz. We will review some of the
classical examples in Section \ref{section3}.

\medskip

For some reason, specialists in the field of minimal surfaces are mostly interested
in embedded surfaces, so the known examples only correspond to
domains with hollow vortices of type I or II.
However, some interesting methods have been developed to
construct minimal surfaces, for example: the {\em conjugate Plateau construction},
 see H. Karcher \cite{karcher}, or the {\em flat structure method} of M. Weber
and M. Wolf \cite{weber-wolf}.
Relaxing the constraint that we want the minimal surface to be embedded, it
should be possible to adapt these methods to construct more general
domains with hollow vortices.
\medskip

In Section \ref{section4}, we discuss another method, which was developed by the author and collaborators to construct minimal surfaces with small catenoidal necks. We will see how it can be
adapted to construct domains with small hollow vortices.
\subsection{Related works}
A family of periodic solutions to the hollow vortex problem was constructed by
Baker, Saffman and Sheffield in \cite{baker}.
It corresponds to the family of horizontal Scherk surfaces
(see Figure \ref{figure2}, top left) -- the authors were of course not aware of that
relationship.
The same solution was derived again by Crowdy and Green in \cite{crowdy},
together with another family of solutions called "staggered vortex streets".
The corresponding minimal graphs are periodic and take on two different values on
the boundary as in Case II.
However, they are asymptotic to half-planes of non-zero slope at
infinity, so extending these surfaces by reflection yields complete minimal surfaces
which are not embedded.

\medskip

Under the additional assumption that $u>0$ in $\Omega$ and $u=0$ on
$\partial\Omega$, solutions to the hollow vortex problem are called
exceptional domains.
Hauswirth, H\'elein and Pacard studied the problem in \cite{hhp}
and discovered an exceptional domain which, as it turns out, corresponds
to the horizontal catenoid.
Partial classification results were obtained by
Khavinson, Lundberg and Teodorescu in \cite{khavinson}.
A complete classification is given in \cite{GAFA}.

\medskip

In a recent paper, Eremenko and Lundberg \cite{eremenko} have investigated
the hollow vortex problem under
the assumption that $u>0$  in $\Omega$ and
$\frac{\partial u}{\partial\nu}>0$ on $\partial\Omega$.
Solutions are called quasi-exceptional domains.
Two examples are constructed. The first one corresponds to
the minimal surface one gets if one tries to add a vertical handle to
the horizontal catenoid.
The second one corresponds to the minimal surface one gets if one tries to
deform the horizontal Scherk surface so that the ``holes'' have different sizes.
Both constructions are known to fail because one cannot solve the vertical period problem.
On the hollow vortex side, this means that the function $u$ takes on different
values on the boundary components, which is perfectly fine.
%
\section{The correspondence}
\label{section2}
\subsection{Preliminary observations}
\label{section21}
We first discuss the hypothesis $||\nabla u||<1$ in the statement of Theorem
\ref{theorem1}.
Let $(\Omega,u)$ be a solution to the hollow vortex problem \eqref{HVP}.
The function $u_z=\frac{1}{2}(u_x-\ii u_y)$ is holomorphic in $\Omega$ and satisfies
$|2u_z|=1$ on $\partial\Omega$.
Let us rule out the trivial case where $u_z$ is constant, in which case $\Omega$ is a half-plane or a band bounded by two parallel lines.
It is known that $||\nabla u||<1$ in $\Omega$ in the following cases:
\begin{enumerate}
\item if $\Omega$ and $u_z$ are doubly periodic, by the maximum principle for holomorphic functions in the quotient,
\item if $||\nabla u||$ is bounded, by a Phragmen Lindel\"of-type result of Fuchs
\cite{fuchs},
\item if $u$ is bounded from below (or above) in $\Omega$, by the proof of
Lemma 2 in \cite{eremenko}.
(In this paper, the authors assume that $\frac{\partial u}{\partial\nu}=+1$ on the boundary, but only the condition $||\nabla u||=1$ is used in the proof of Lemma 2.)
\end{enumerate}
Let us also mention that provided $||\nabla u||<1$, the domain $\Omega$ must be strictly concave: see the proof of Proposition 4 in \cite{GAFA}.
\subsection{Weierstrass representation}
For the reader not familiar with minimal surfaces and to fix notations, we
recall the Weierstrass representation formula: \begin{equation}
\label{weierstrass}
X(z)=(x_1(z),x_2(z),x_3(z))=X_0+\Re\int_{z_0}^z \left(\frac{1}{2}(g^{-1}-g)\omega,\frac{\ii}{2}(g^{-1}+g)\omega,\omega\right)
\end{equation}
In its local form, $g$ is a meromorphic function on a simply connected domain
$\Sigma\subset\C$ and $\omega=f(z)dz$ where $f$ is a holomorphic function on
$\Sigma$ having a zero at each zero or pole of $g$, with the same multiplicity.
$z_0\in\Sigma$ is an arbitrary base point and $X_0$ is some constant vector.
Then $X:\Sigma\to \R^3$ is a conformal parametrization of a minimal surface $M$.
Moreover, the Gauss map of $M$ is given by
$$N=\left(\frac{2\,\Re (g)}{|g|^2+1},\frac{2\,\Im (g)}{|g|^2+1},\frac{|g|^2-1}{|g|^2+1}\right).$$
In other words, $g$ is the stereographic projection of the Gauss map.
\medskip

In its global form, $\Sigma$ is a Riemann surface, $g$ is a meromorphic function and $\omega$ is a holomorphic 1-form on $\Sigma$.
\subsection{The correspondence "vortex $\to$ minimal"}
Let $(\Omega,u)$ be a solution to the hollow vortex problem \eqref{HVP}.
Consider the minimal surface $M$ given by the Weierstrass
representation formula \eqref{weierstrass} with
$$g=\displaystyle\frac{-1}{2u_z},\quad\omega=2u_z\,dz,\quad
X_0=(0,0,u(z_0)).$$
Let $\psi(z)=x_1(z)+\ii x_2(z):\Omega\to\C$.
\begin{proposition}
\label{correspondence1}
In the above setup:
\begin{enumerate}
\item $x_3(z)=u(z)$.
\item $\psi(z)$ is well defined in $\Omega$, namely does not depend on the integration path
from $z_0$ to $z$.
\item $d\psi=dz$ along $\partial\Omega$.
\end{enumerate}
Assume moreover that $||\nabla u||<1$ in $\Omega$. Then
\begin{enumerate}
\setcounter{enumi}{3}
\item For any $z'\neq z$ in $\Omega$, $0<|\psi(z')-\psi(z)|<|z'-z|$.
\item $\psi$ is a diffeomorphism from $\Omega$ to $\whOmega=\psi(\Omega)$.
\item The boundary of $\whOmega$ is $\psi(\partial\Omega)$.
\end{enumerate}
\end{proposition}
The proof of this proposition is essentially the same as the proof of Theorem 9
in \cite{GAFA}. For completeness, we give the details in Appendix \ref{appendix1}.
\medskip

From Points (1) and (5), we see that $M$ is the graph of
$v=u\circ\psi^{-1}$ over the domain $\whOmega$.
Since $|g|=1$ on $\partial\Omega$, the Gauss map is horizontal on the boundary, so
$||\nabla v||\to\infty$ on $\partial\whOmega$.
From Point (2), we see that each component of $\partial\whOmega$ is a translation
of the corresponding component of $\partial\Omega$. Moreover, from Point (4),
we see that $\psi$ moves the boundary components toward each other.
\begin{remark} In \cite{GAFA} we took $g=2u_z$, so $|g|<1$ in $\Omega$
and the Gauss map of $M$ was pointing down. For a graph it is more natural to
choose the upward pointing normal so that the horizontal projection preserves 
orientation. The correspondence has better properties with the antipodal  choice
$g=\frac{-1}{2u_z}$.
\end{remark}
\subsection{The correspondence "minimal $\to$ vortex"}
\label{minimal-vortex}
Let $(\whOmega,v)$ be a solution to Problem \eqref{mse} and $M$ be the
minimal surface given as the graph of $v$. We orient $M$ by its upward pointing
normal. Then $M$ is parametrized on some Riemann surface (with boundary)
$\Sigma$ by the Weierstrass representation formula \eqref{weierstrass}.
We have $|g|>1$ on $\Sigma$ and $|g|=1$ on $\partial\Sigma$.
We observe that even though $\Sigma$ is diffeomorphic to the planar domain
$\whOmega$, in practice, it will not be given explicitely as a domain in the plane
(see examples in Section \ref{section3}), so it is better to leave it as an abstract
Riemann surface.
Let $\psi(z)=x_1(z)+\ii x_2(z)$. Since $M$ is a graph, $\psi$ is a diffeomorphism
from $\Sigma$ to $\whOmega$. Define $F:\whOmega\to\C$ by
$$F(\psi(z))=-\displaystyle\int_{z_0}^z g\omega.$$
\begin{proposition}
\label{correspondence2}
In the above setup:
\begin{enumerate}
\item $F$ is well defined in $\whOmega$.
\item $dF=dz$ on $\partial\whOmega$.
\item For any $z\neq z'$ in $\whOmega$, $|F(z)-F(z')|> |z-z'|$.
\item $F$ is a diffeomorphism from $\whOmega$ to $\Omega=F(\whOmega)$.
\item The function $u(z)=v(F^{-1}(z))$ solves Problem \eqref{HVP}
and satisfies $||\nabla u||<1$ in $\Omega$.
\end{enumerate}
\end{proposition}
The proof of this proposition is essentially the same as the proof of Theorem 10
in \cite{GAFA}. For completeness, we give the details in Appendix \ref{appendix2}.
\medskip

The maps
$(\Omega,u)\mapsto (\whOmega,v)$
and $(\whOmega,v)\mapsto(\Omega,u)$ defined by Propositions \ref{correspondence1}
and \ref{correspondence2} are inverse of each other,
provided we identify two domains which
differ by a translation. See Theorem 11 in \cite{GAFA}.
\begin{remark} A computation shows that $dF$ is given in term of $v$ by
$$dF=dx+\ii dy+\frac{(1+v_x^2)dx+v_x v_ydy}{W}+\ii\,\frac{ v_x v_y dx+(1+v_y^2)dy}{W}$$
where $W=\sqrt{1+v_x^2+v_y^2}.$
Alternately, one could take this as a definition of $F$. (The minimal surface equation
implies that $dF$ is closed.)
With this definition, the proof of
Proposition \ref{correspondence2} is more computational but
avoids Weierstrass representation.
\end{remark}
\section{Classical examples}
\label{section3}
\begin{figure}
\begin{center}
\includegraphics[height=60mm,angle=90]{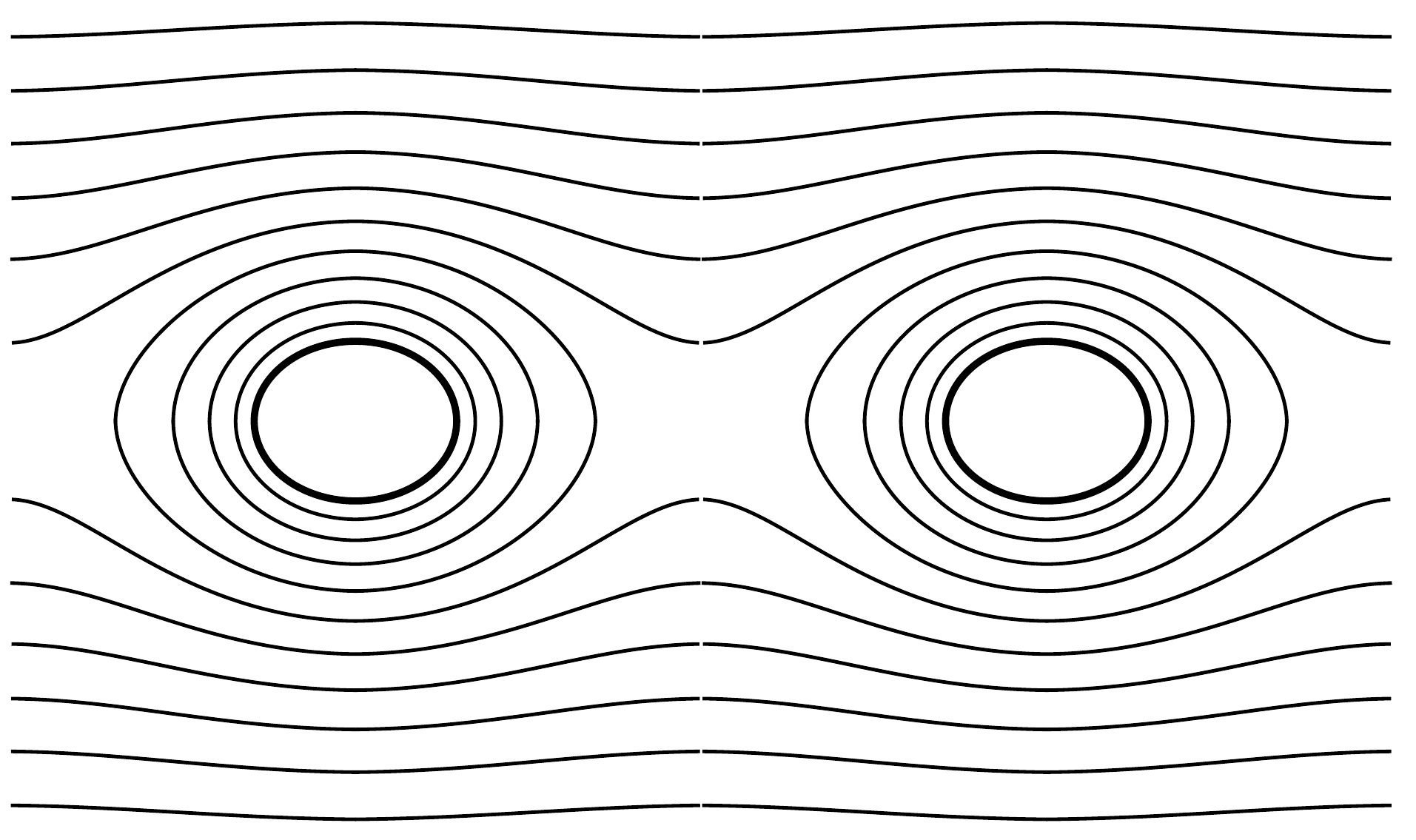}
\hspace{1.9cm}
\includegraphics[height=40mm,angle=90]{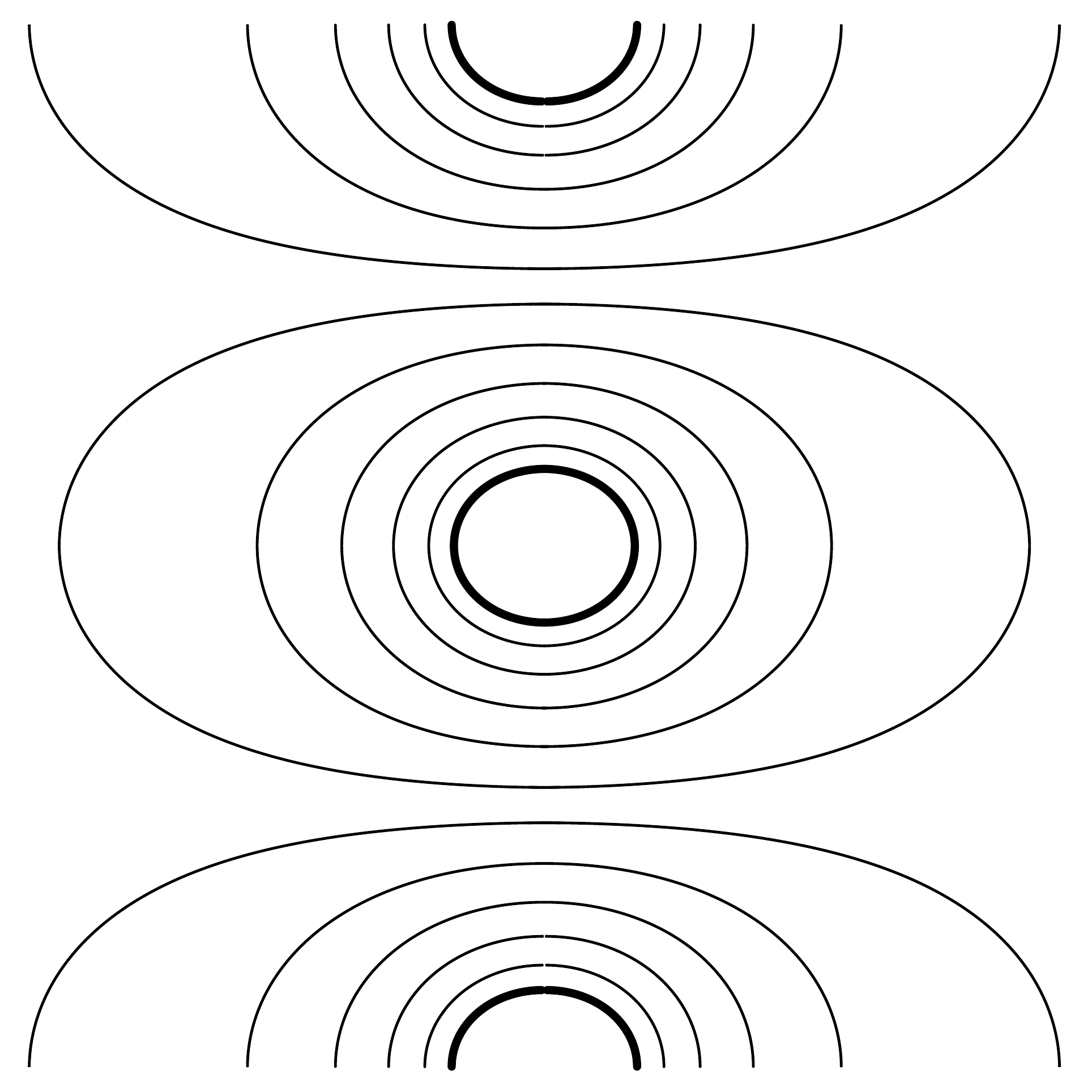}
\vspace{2.0cm}

\includegraphics[height=50mm]{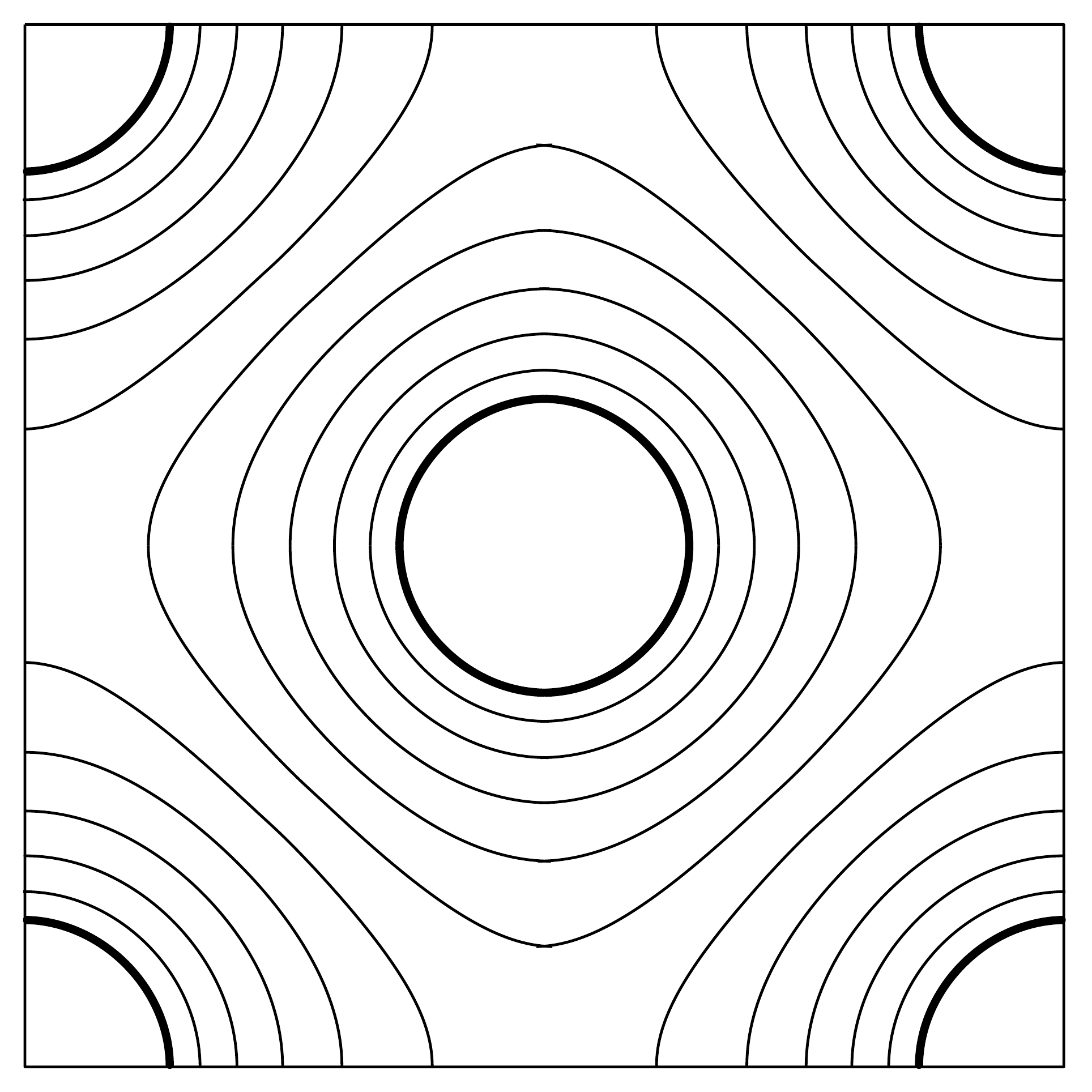}
\hspace{1.9cm}
\includegraphics[height=50mm]{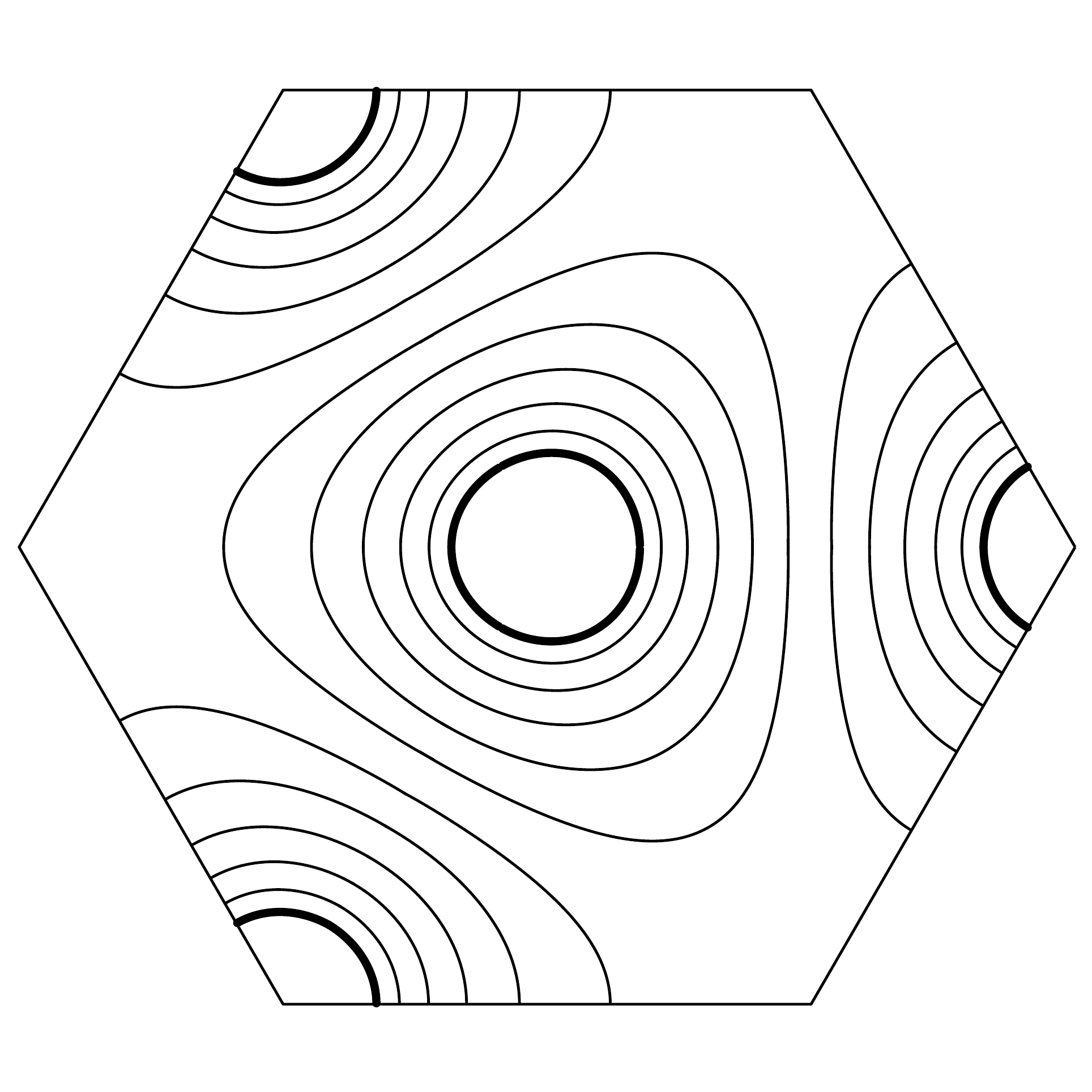}
\end{center}
\caption{Top: the periodic domains corresponding
to the horizontal Scherk surface (left) and a toroidal half-plane layer with $a=0.5$ (right). 
Bottom: the doubly-periodic domains corresponding to Schwarz P-surface
(left) and H-surface with $a=0.5$ (right). The curves are the stream lines. Computed with Maple.}
\label{figure2}
\end{figure}

We focus on examples which are bounded by closed curves, as they are probably more interesting from the hydrodynamics point of view.
All these examples admit deformations, and a lot more examples are known, see
\cite{karcher}.
The Weierstrass data for all these examples is explicit and one can compute
numerically the corresponding domain with hollow vortices: see Figure \ref{figure2}.

\begin{enumerate}
\item The horizontal Scherk surface, a periodic minimal surface with horizontal period:
$$ g=\frac{1}{z},\quad \omega=\frac{z\,dz}{z^4+6z^2+1}.$$
\item Karcher toroidal halfplane layers, a family of doubly periodic minimal surfaces with one horizontal and one vertical periods:
$$g=\frac{1}{z},\quad \omega=\frac{dz}{\sqrt{(z^2+a^2)(z^2+a^{-2})}},\quad 0<a<1.$$
\item Schwarz P-surface, a triply periodic minimal surface:
$$g=\frac{1}{z},\quad \omega=\frac{z\,dz}{\sqrt{z^8-14 z^4+1}}.$$
\item Schwarz H-surfaces, a family of triply periodic minimal surfaces:
$$g=\frac{1}{z},\quad \omega=\frac{z\,dz}{\sqrt{z(z^3+a^3)(z^3+a^{-3})}},\quad 0<a<1.$$
\end{enumerate}
In all these examples, the Riemann surface $\Sigma$ is a branched cover of the unit
disk $D(0,1)$ punctured at $z=\pm \ii(\sqrt{2}-1)$ in Case (1) and
$z=0$ in Case (2). The residues at the punctures and the multivaluation of the
square roots are responsible for the periods of the minimal surfaces and the
corresponding domains $\Omega$.

\section{Domains with small hollow vortices}
In this section, we give a general method to construct domains with small hollow
vortices, first in the finite connectivity case, then in the periodic case.
In what follows, we identify points and vectors in the plane $\R^2$ with complex
numbers.
\label{section4}
\subsection{Domains with a finite number of hollow vortices}
First some definitions.
Let $(\Omega,u)$ be a solution to Problem \eqref{HVP}.
Let $\gamma$ be a closed curve in the boundary of $\Omega$.
Let  $\boC(\gamma)=\int_{\gamma} \vec{v}\cdot \vec{d\ell}$ be the circulation of $\vec{v}$ on the curve $\gamma$.
Since $||\vec{v}||=1$ on $\gamma$, the absolute value of $\boC(\gamma)$ is equal to the length of $\gamma$, but it can have either sign, depending on wether
the vortex is ``spinning'' left or right.
\begin{definition}
A vortex configuration is a finite set of $n\geq 2$ distinct points $p_1,\cdots,p_n$ in the complex plane with weights $c_1,\cdots,c_n$ which are non-zero real numbers.
Forces are defined by
$$F_i=\sum_{j\neq i} \frac{c_ic_j}{p_i-p_j}.$$
We say a configuration is balanced if $F_i=0$ for $i=1,\cdots,n$.
We say a balanced configuration is non-degenerate if the $n\times n$ jacobian matrix
$\frac{\partial F_i}{\partial p_j}$ has complex rank $n-2$.
\end{definition}
Observe that we always have
\begin{equation}
\label{sum-Fi}
\sum_{i=1}^n F_i=0
\end{equation}
\begin{equation}
\label{sum-piFi}
\sum_{i=1}^n p_i F_i=\sum_{i<j} c_i c_j.
\end{equation}
Hence $n-2$ is the maximum rank that the jacobian matrix may have.
Also, \eqref{sum-piFi} gives a restriction on the weights for a balanced configuration
to exist.
\begin{theorem}
\label{theorem-finite}
Given a balanced, non-degenerate configuration, there exists a 1-parameter
family of solutions $(\Omega_t,u_t)$ of the hollow vortex problem \eqref{HVP},
depending on a small parameter $t>0$, such that:
\begin{enumerate}
\item $\Omega_t$ has $n$ boundary components, denoted 
$\gamma_{1,t}\cdots,\gamma_{n,t}$, all of them closed curves.
\item The circulation $\boC(\gamma_{i,t})$ is equal to $2\pi c_i t$.
\item As $t\to 0$, $\gamma_{i,t}$ shrinks to the point $p_i$. Moreover, its
asymptotic shape is circular.
\end{enumerate}
\end{theorem}
Here is a simple example of balanced configuration, with
dihedral symmetry of order $n-1$ ($n\geq 3$):
$$c_j=1,\quad p_j=e^{2\pi \ii j/(n-1)}\quad \mbox{ for $1\leq j\leq n-1$},$$
$$c_n=1-\frac{n}{2},\quad p_n=0.$$
The configuration is balanced by symmetry and Equation \eqref{sum-piFi}.
One can easily break the symmetries by perturbing the weights.
\medskip

The proof of Theorem \ref{theorem-finite} follows \cite{nosym} in the minimal case
and is omitted. It is also very similar to the proof of Theorem \ref{theorem-periodic}
in Appendix \ref{appendix3}.
\subsection{Periodic domains}
In this section, we construct periodic domains with hollow vortices
and period $T$ in the $x$-direction,
and such that the velocity vector has a limit as $y\to\pm\infty$.
The limit velocities cannot be arbitrary, as the following proposition shows: 
\begin{proposition}
\label{proposition-ab}
Let $\Omega$ be a domain with hollow vortices. Assume that
\begin{enumerate}
\item $\Omega$ and the velocity vector $\vec{v}$ are periodic with period $T$ in the
$x$-direction.
\item The quotient $\Omega/T$ has a finite number of boundary components,
denoted $\gamma_1,\cdots,\gamma_n$, all of them closed curves.
\item The velocity vector $\vec{v}$ has a limit as $y\to+\infty$ and $y\to-\infty$, denoted respectively
$v^+$ and $v^-$.
\end{enumerate}
Then either:
\begin{itemize}
\item[(a)] $v^+=v^-$ and $\displaystyle\sum_{i=1}^n\boC(\gamma_i)=0$, or
\item[(b)] $v^+=-v^-=\displaystyle\frac{1}{2T}\sum_{i=1}^n \boC(\gamma_i)$.
\end{itemize}
\end{proposition}
Here we see the vectors $v^+$ and $v^-$ as complex numbers. In Case (b), the
limit velocity must be real numbers, while there is no such restriction in
Case (a).

\medskip
Proof. Since we identify vectors with complex numbers, we can write
$\vec{v}=-2\ii u_{\overline{z}}$. Hence
$$\lim_{y\to\pm\infty}2u_z=-\ii\overline{v^{\pm}}.$$
If $\gamma$ is a stream line, the circulation of the velocity is related to $u_z$ by
\begin{equation}
\label{circulation}
\int_{\gamma} 2u_z\,dz=\int_{\gamma} u_x\,dx+u_y\,dy
-\ii\int_{\gamma} u_y\,dx-u_x\,dy
=-\ii\,\boC(\gamma).
\end{equation}

For large $R$, consider the domain $\Omega_R=\Omega\cap\{ -R<y<R\}$. By Cauchy Theorem,
$$0=\int_{\partial(\Omega_R/T)}2u_z\,dz=
\sum_{i=1}^n\int_{\gamma_i}2u_z\,dz+\int_{z=-\ii R}^{T-\ii R} 2u_z\,dz
+\int_{z=T+\ii R}^{\ii R} 2u_z\,dz.$$
We let $R\to\infty$ and obtain
\begin{equation}
\label{cauchy1}
\sum_{i=1}^n\boC(\gamma_i)=T(\overline{v^+}-\overline{v^-}).
\end{equation}
Since $u$ is constant on $\gamma_i$, we have $du=u_z\,dz+u_{\overline{z}}\,d\overline{z}=0$ along $\gamma_i$. Hence using $|2u_z|=1$ on $\gamma_i$,
$$\int_{\gamma_i} (2u_z)^2dz=-\int_{\gamma_i} 4u_z\overline{u_z}\,d\overline{z}=-\int_{\gamma_i}
d\overline{z}=0.$$
Using Cauchy theorem again with the function $(2u_z)^2$ and letting
$R\to\infty$ gives
\begin{equation}
\label{cauchy2}
T((\overline{v^+})^2-(\overline{v^-})^2)=0.
\end{equation}
Proposition \ref{proposition-ab} follows from \eqref{cauchy1} and \eqref{cauchy2}.
\cqfd
\begin{definition}
\label{definition-periodic-configuration}
A periodic vortex configuration is a finite set of non-zero complex numbers $p_1,\cdots,
p_n$ with weight $c_1,\cdots,c_n$ which are non-zero real numbers, together with a non-zero complex number $c_0$.
We assume that either:
\begin{itemize}
\item[]Case (a) $c_1+\cdots+c_n=0$, or
\item[]Case (b) $c_1+\cdots+c_n+2c_0=0$.
\end{itemize}
We define forces by
$$F_i=\sum_{j\neq  i}c_i c_j \frac{p_i+p_j}{p_i-p_j}+\left\{\begin{array}{ll}
2c_ic_0 & \mbox{ in Case (a)}\\
0 &\mbox{ in Case (b)}
\end{array}\right.$$
We say a configuration is balanced if $F_i=0$ for $1\leq i\leq n$. 
We say a balanced configuration is non-degenerate if the jacobian matrix 
$\frac{\partial F_i}{\partial p_j}$ has complex rank $n-1$.
\end{definition}
Observe that in either case, $F_1+\cdots+F_n=0$, so $n-1$ is the maximum rank
that the jacobian matrix may have.
\begin{theorem}
\label{theorem-periodic}
Given a balanced, non-degenerate periodic configuration, there exists a 1-parameter
family of solutions $(\Omega_t,u_t)$ of the hollow vortex problem \eqref{HVP},
depending on a small parameter $t>0$, such that:
\begin{enumerate}
\item $\Omega_t$ is a periodic domain with period $T=2\pi$.
\item The quotient $\Omega_t/T$ has $n$ boundary components, denoted 
$\gamma_{1,t}\cdots,\gamma_{n,t}$, all of them closed curves.
\item The circulation $\boC(\gamma_{i,t})$ is equal to $2\pi c_i t$.
\item As $t\to 0$, $\gamma_{j,t}$ shrinks to the point $q_j=\ii\log p_j$. Moreover, its
asymptotic shape is circular.
\item In Case (a), the limit of the velocity as $y\to\pm\infty$ is $t\overline{c_0}$.
In Case (b), the limit of the velocity as $y\to\pm\infty$ is $\mp tc_0$.
\end{enumerate}
\end{theorem}
We prove this theorem in Appendix \ref{appendix3}. The proof follows \cite{connor}
in the minimal case.
Please take care that the limit position of the vortices is $q_j=\ii\log p_j$ and
not $p_j$ as in Theorem \ref{theorem-finite}. The $2\pi \ii$ multivaluation of the complex
logarithm is responsible for the period $T=2\pi$ of the domain.
In term of the points $q_j$, the forces are given by
$$F_i=-\ii\sum_{j\neq  i}c_i c_j \cot\frac{q_i-q_j}{2}+\left\{\begin{array}{ll}
2c_ic_0 & \mbox{ in Case (a)}\\
0 &\mbox{ in Case (b)}
\end{array}\right.$$
\subsection{Examples of periodic configurations of type (a)}
We take $n=2$, $c_1=1$ and $c_2=-1$. Solving $F_1=0$ gives
$$q_2=q_1+\ii\log \frac{2c_0-1}{2c_0+1},\qquad c_0\neq \pm 0.5.$$
The points $q_i$ for various values of $c_0$
are represented on Figures \ref{fig-a-1}, \ref{fig-a-2}, \ref{fig-a-3}.

\medskip

If we take $n=5$, $c_1=c_2=c_3=1$ and $c_4=c_5=-1.5$, we obtain
an uneven vortex street: see Figure \ref{fig-a-4}.
\subsection{Examples of periodic configurations of type (b)}
First assume that all $c_i$ are equal to $1$. The configuration
$p_j=e^{2\pi \ii j/n}$ is balanced by symmetry: all forces $F_i$ are equal and their
sum is zero.
This gives $q_j=-2\pi j/n$ so the vortices are regularly spaced. This configuration
gives the family of domains corresponding to the family of horizontal Scherk surfaces
when it is close to its catenoidal limit.

\medskip

To get a more interesting example, take $n=3$, $c_1=c_2=1$ and leave
$c_3$ as a parameter. We may
normalize $p_3=1$. Computations show that $p_1$ and $p_2$ are the roots
of the polynomial
$\displaystyle P(z)=z^2+\frac{2c_3}{c_3+1}z+1.$
We obtain a two lanes vortex street, see Figure \ref{fig-b}.

\begin{figure}
\begin{center}
\setlength{\unitlength}{2.0cm}
\begin{picture}(6,0.8)(0,-0.5)
\put(0,0){$\circ$}
\put(1,-0.3497){$\bullet$}
\put(2,0){$\circ$}
\put(3,-0.3497){$\bullet$}
\put(4,0){$\circ$}
\put(5,-0.3497){$\bullet$}
\put(3.5,0.2){$\longrightarrow$}
\put(3.5,-0.5){$\longrightarrow$}
\end{picture}
\end{center}
\caption{A periodic configuration of type (a) with $n=2$, $c_0=0.25$. Circles represent
vortices with right spin ($c_i>0$), bullets represent vortices with left spin
($c_i<0$). The arrows indicate the direction of the velocity at infinity.
Three fundamental domains are represented. Compare with Figure \ref{fig1}.}
\label{fig-a-1}

\begin{center}
\setlength{\unitlength}{2.0cm}
\begin{picture}(6,0.8)(0,-0.5)
\put(0,0){$\circ$}
\put(0,-0.3497){$\bullet$}
\put(2,0){$\circ$}
\put(2,-0.3497){$\bullet$}
\put(4,0){$\circ$}
\put(4,-0.3497){$\bullet$}
\put(3.5,0.2){$\longrightarrow$}
\put(3.5,-0.5){$\longrightarrow$}
\end{picture}
\end{center}
\caption{A periodic configuration of type (a) with $n=2$, $c_0=1$.}
\label{fig-a-2}

\begin{center}
\setlength{\unitlength}{2.0cm}
\begin{picture}(6,0.8)(0,-0.2)
\put(0,0){$\circ$}
\put(0.2406,-0.2041){$\bullet$}
\put(2,0){$\circ$}
\put(2.2406,-0.2041){$\bullet$}
\put(4,0){$\circ$}
\put(4.2406,-0.2041){$\bullet$}
\put(1.5,0.2){$\nearrow$}
\put(3.5,0.2){$\nearrow$}
\end{picture}
\end{center}
\caption{A periodic configuration of type (a) with $n=2$, $c_0=e^{-\ii\pi/4}$.}
\label{fig-a-3}

\begin{center}
\setlength{\unitlength}{2.0cm}
\begin{picture}(6,0.8)(-1,-0.5)
\put(-.6666666664, -0.02936340626){$\circ$}
\put(-0.03551828126, 0.01468170323){$\circ$}
\put(.7021849478, 0.01468170323){$\circ$}
\put(1.333333334, -0.02936340626){$\circ$}
\put(1.964481719, 0.01468170323){$\circ$}
\put(2.702184948, 0.01468170323){$\circ$}
\put(3.333333334, -0.02936340626){$\circ$}
\put(3.964481719, 0.01468170323){$\circ$}
\put(4.702184948, 0.01468170323){$\circ$}
\put(-.1846663713, -.5251421936){$\bullet$}
\put(.8513330378, -.5251421936){$\bullet$}
\put(1.815333629, -.5251421936){$\bullet$}
\put(2.851333038, -.5251421936){$\bullet$}
\put(3.815333629, -.5251421936){$\bullet$}
\put(4.851333038, -.5251421936){$\bullet$}
\put(2.5,0.2){$\longrightarrow$}
\put(2.5,-0.7){$\longrightarrow$}
\end{picture}
\end{center}
\caption{A periodic configuration of type (a) with $n=5$, $c_1=c_2=c_3=1$, $c_4=c_5=-1.5$ and $c_0=0.5$. Computed with Maple.}
\label{fig-a-4}

\begin{center}
\setlength{\unitlength}{1.5cm}
\begin{picture}(6,2)(-1,-1)
\put(-.9999999999, -.5610998522){$\circ$}
\put(-.9999999999, .5610998522){$\circ$}
\put(1.000000000, -.5610998522){$\circ$}
\put(1.000000000, .5610998522){$\circ$}
\put(3.000000000, -.5610998522){$\circ$}
\put(3.000000000, .5610998522){$\circ$}
\put(0,0){$\bullet$}
\put(2,0){$\bullet$}
\put(4,0){$\bullet$}
\put(2.5,0.8){$\longrightarrow$}
\put(2.5,-0.8){$\longleftarrow$}
\end{picture}
\end{center}
\caption{A periodic configuration of type (b) with $n=3$, $c_1=c_2=1$ and $c_3=-1.5$.}
\label{fig-b}
\end{figure}

\appendix
\section{Proof of Proposition \ref{correspondence1}}
\label{appendix1}
Proof of Point (1):
$$x_3(z)=u(z_0)+\Re\int_{z_0}^z 2u_z\,dz
=u(z_0)+\int_{z_0}^z (u_z\,dz+u_{\overline{z}}\,d\overline{z})
=u(z_0)+\int_{z_0}^z du
=u(z).$$
Proof of Point (2): consider the differential
$$d\psi=dx_1+\ii dx_2=\Re\left(\frac{1}{2}(g^{-1}-g)\omega\right)+\ii\,\Re\left(\frac{\ii}{2}(g^{-1}+g)\omega\right)
=\frac{1}{2}(\overline{g^{-1}\omega}-g\omega).$$
With our choice of $g$ and $\omega$,
\begin{equation}
\label{dpsi}
d\psi=\frac{1}{2}\left(dz-4(u_{\overline{z}})^2\,d\overline{z}\right).
\end{equation}
We have to prove that $d\psi$ is an exact differential (namely, the differential of
a globally defined function $\psi$).
If $t\mapsto\gamma(t)$ is a parametrization of a boundary component of $\Omega$,
then since $u$ is constant on $\gamma$:
$$du(\gamma')=0=\left(u_z dz+u_{\overline{z}}d\overline{z}\right)(\gamma').$$
\begin{equation}
\label{dpsi2}
d\psi(\gamma')=\frac{1}{2}\left(dz(\gamma')+4u_{\overline{z}}u_z\,dz(\gamma')\right)
=\frac{1}{2}(1+||\nabla u||^2)dz(\gamma')=dz(\gamma').
\end{equation}
Hence if $\gamma$ is a closed component of $\partial\Omega$,
$\int_{\gamma}d\psi=0$. Since $\Omega$ is a planar domain, this implies that
$d\psi$ is an exact differential.
Also \eqref{dpsi2} proves Point (3).
\medskip

Proof of Point (4): We prove that
\begin{equation}
\label{zz}
|2(\psi(z')-\psi(z))-(z'-z)|<|z'-z|
\end{equation}
which implies Point (4) by triangular inequality.
We may decompose the segment $[z,z']$ into $n$ segments
$[z_i,z_{i+1}]$ for $1\leq i\leq n$, such that $z_1=z$, $z_{n+1}=z'$,
$z_i\in\partial\Omega$ for $2\leq i\leq n$
and for each $i$, the open segment $(z_i,z_{i+1})$ is either included in
$\Omega$ or its complement.
In the first case, we have by Equation \eqref{dpsi} and using $||\nabla u||<1$
\begin{equation}
\label{zz2}
|2(\psi(z_{i+1})-\psi(z_i))-(z_{i+1}-z_i)|
=\left|\int_{z_i}^{z_{i+1}} 4(u_{\overline{z}})^2\,d\overline{z}\right|
\leq \int_{z_i}^{z_{i+1}}4|u_{\overline{z}}|^2\,|dz|<|z_{i+1}-z_i|.
\end{equation}
In the second case, since $\Omega$ is a concave domain, $z_i$ and $z_{i+1}$ must
be on the same component of $\partial\Omega$. By Point (3),
$\psi(z_{i+1})-\psi(z_i)=z_{i+1}-z_{i}$ so \eqref{zz2} becomes an equality.
Summing for $1\leq i\leq n$ gives \eqref{zz}.
\medskip

Proof of Point (5): Equation \eqref{dpsi} and $||\nabla u||<1$ implies that
$d\psi$ is an isomorphism so $\psi$ is a local diffeomorphism.
Point (4) implies that $\psi$ is injective, so is a global diffeomorphism
onto its image.
\medskip

Proof of Point (6): since $\psi$ is a homeomorphism from $\Omega$ to
$\whOmega$ and extends continuously to $\overline{\Omega}$,
we have $\psi(\partial\Omega)\subset\partial\whOmega$ by elementary
topology. (Here $\overline{\Omega}$ denotes the closure of $\Omega$.)
Assume by contradiction that $\psi(\partial\Omega)\neq\partial\whOmega$
and let $a_0\in \partial\whOmega\setminus\psi(\partial\Omega)$.
\medskip

The finiteness hypothesis (Hypothesis \ref{finiteness}) ensures that
$\psi(\partial\Omega)$ is closed. Indeed, for each component $\gamma$
of $\partial\Omega$, $\psi(\gamma)$ is a translate of $\gamma$ so
is closed, and the finiteness hypothesis prevent them from accumulating.
\medskip

Let $\varepsilon=d(a_0,\psi(\partial\Omega)$. Choose a point $a_1\in\whOmega$
such that $|a_0-a_1|\leq \frac{\varepsilon}{4}$.
Let $a_2$ be a point on $\partial\whOmega$ whose distance to $a_1$ is minimum.
Then $d(a_2,\psi(\partial\Omega))\geq\frac{\varepsilon}{2}$ and the semi-open
segment $[a_1,a_2)$ is entirely included in $\whOmega$.
Let $\alpha(t):[0,\ell)\to\Omega$ be a path such that $\psi(\alpha(t))$ is the parametrization
at unit speed of the segment $[a_1,a_2)$.
We must have $||\alpha(t)||\to\infty$ as $t\to\infty$, else $a_2$ would be in
$\psi(\overline{\Omega})$.
This implies that the path $\alpha$ has infinite length.
Now the conformal metric induced by the minimal immersion $X$ is given by
$$ds=\frac{1}{2}(|g|+|g|^{-1})|\omega|\geq \frac{1}{2}|dz|.$$
Hence, the curve $X(\alpha(t))$ on $M$ has infinite length. This curve is the graph of
the function $v$ on the segment $[a_1,a_2)$.
By standard results in minimal surface theory (see the proof of Lemma 2
in \cite{GAFA} for the details), 
this implies that $\lim_{z\to a_2} v(z)=\pm\infty$.
Moreover, there exists a divergence line $L$, containing $a_2$ and contained in $\partial\whOmega$, such that $v\to\pm\infty$ on $L$.
By connectedness, $\whOmega$
must be on one side of $L$. To derive a contradiction, we distinguish two cases:
\begin{itemize}
\item If $d(L,\psi(\partial\Omega))>0$, then the function $v$ satisfies the minimal surface
equation in a band with boundary value $\pm\infty$ on one side. This is impossible
by Proposition 1 in \cite{mazet}.
\begin{remark} If all components of $\partial\Omega$ are closed curves, then the finiteness hypothesis implies that $d(L,\psi(\partial\Omega))>0$.
\end{remark}
\item If $d(L,\psi(\partial\Omega)=0$, then
there is an unbounded component of $\partial\Omega$, say $\gamma_1$, such that
$\psi(\gamma_1)$ is asymptotic to $L$. Also, there can be at most two such components.
Label $\gamma_2$ the other component asymptotic to $L$, if any.
Let us write $\whgamma_i=\psi(\gamma_i)$.
There exists $\varepsilon>0$ so that all other components of $\psi(\partial\Omega)$
are at distance greater than $\varepsilon$ of $L$. We obtain a contradiction using
the catenoid as a barrier as in the proof of the strong half-space theorem of
Hoffman Meeks \cite{hoffman-meeks}. The only difference is that $M$ is not complete
so we have to mind its boundary.

 \medskip

Without loss of generality, we may assume that $L$ is the line $x_2=0$ in the horizontal
plane, $M$ lies in the half-space $x_2>0$, and also $v<0$ on $\whgamma_1$ and $\whgamma_2$ and
$v\to+\infty$ on $L$.
Let $\nu$ be the interior conormal to the boundary of $M$.
Since $\whgamma_1$ and $\whgamma_2$ are horizontal symmetry curves, $\nu$
is vertical. Evaluating the vertical flux in the subdomain of $\whOmega$
defined by $R<|x_1|<2R$ and $0<x_2<\varepsilon$ for large values of $R$, we obtain that $\nu=(0,0,1)$ on $\whgamma_1$ and $\whgamma_2$. 
\medskip

Let $C$ be the horizontal half-catenoid $x_1^2+x_3^2=\cosh^2 x_2$, $x_2<0$.
Let $C_t=(0,a,0)+t\,C_1$, where $0<a<\varepsilon$ and $0<t\leq 1$.
If $a$ is small enough then $C_1$ does not intersect $M$. Also, as $t\to 0$,
$C_t$ converges to the vertical plane $x_2=a$, so $C_t$ intersects $M$ for $t>0$
small enough. Let $t_0<1$ be the largest time so that $C_t$ intersects $M$.
Then $C_{t_0}$ intersects $M$ at a boundary point. Since $C_t$ lies in the half-space
$x_2<a<\varepsilon$, that point must be on $\whgamma_1$ or $\whgamma_2$.
Since $x_3<0$ on $\whgamma_i$ and $\nu=(0,0,1)$, $M$ will still intersect $C_t$
for $t$ slightly larger than $t_0$, a contradiction.\cqfd
\end{itemize}
\section{Proof of Proposition \ref{correspondence2}}
\label{appendix2}
Proof of Point (1): consider the differential $d\varphi=-g\omega$ on $\Sigma$.
We have to prove that $d\varphi$ is an exact differential.
Since $\omega$ has a zero at each pole of $g$, $d\varphi$ is holomorphic
in $\Sigma$.
Let $\gamma$ be a component of $\partial\Omega$. Since $|g|=1$ on
$\gamma$ and $\omega(\gamma')$ is imaginary,
\begin{equation}
\label{dphi}
d\psi(\gamma')=\frac{1}{2}\left(\overline{g^{-1}\omega(\gamma')}-g\omega(\gamma')\right)
=-g\omega(\gamma')=d\varphi(\gamma').
\end{equation}
Since $d\psi$ is an exact differential and $\Sigma$ is diffeomorphic to a planar
domain, $d\varphi$ is the differential of a globally defined function $\varphi$.
Then $F=\varphi\circ\psi^{-1}$ is well defined.
Point (2) is a consequence of \eqref{dphi}.
\medskip

Proof of Point (3): let $\tau$ be the unit vector in the direction of $z'-z$.
We prove that
\begin{equation}
\label{Ftau}
\langle F(z')-F(z),\tau\rangle>\langle z'-z,\tau\rangle
\end{equation}
which implies Point (3).
We may decompose the segment $[z,z']$ into $n$ segments
$[z_i,z_{i+1}]$ for $1\leq i\leq n$, such that $z_1=z$, $z_{n+1}=z'$,
$z_i\in\partial\whOmega$ for $2\leq i\leq n$
and for each $i$, the open segment $(z_i,z_{i+1})$ is either included in
$\whOmega$ or its complement.
In the first case, let $\alpha(t):(0,\ell)\to\Sigma$ be such that
$\psi(\alpha(t))$ is the parametrization of the segment $(z_i,z_{i+1})$
at constant speed $\tau$, in other words $d\psi(\alpha')=\tau$.
Then
\begin{eqnarray*}
\langle d\varphi(\alpha'),\tau\rangle
&=&\langle d\varphi(\alpha'),d\psi(\alpha')\rangle\\
&=&\langle d\varphi(\alpha')-d\psi(\alpha'),d\psi(\alpha')\rangle +||d\psi(\alpha')||^2\\
&=& \frac{1}{4}\left(|g\omega(\alpha')|^2-|g^{-1}\omega(\alpha')|^2\right)+1\\
&>&1 \quad\mbox{ since $|g|>1$}.
\end{eqnarray*}
Integrating from $t=0$ to $\ell$, we obtain
\begin{equation}
\label{Ftau2}
\langle F(z_{i+1})-F(z_i),\tau\rangle > \ell =\langle z_{i+1}-z_i,\tau\rangle.
\end{equation}
If the segment $(z_i,z_{i+1})$ is included in the complementary of $\whOmega$,
then since $\whOmega$ is a concave domain, $z_i$ and $z_{i+1}$ must be
on the same component of $\partial\whOmega$, so $F(z_i)=F(z_{i+1})$ by
Point (2).
Hence \eqref{Ftau2} become an equality.
Summing for $1\leq i\leq n$ gives \eqref{Ftau}.
\medskip

Proof of Point (4): since $|g|>1$ in $\Sigma$, $d\varphi\neq 0$ in $\Sigma$ so
$\varphi$ and $F$  are local diffeomorphisms. By Point (3), $F$ is injective, so
is a global diffeomorphism onto its image $\Omega$. By Point (3),
$F$ is proper, so $\partial\Omega=F(\partial\whOmega)$.
\medskip

Proof of Point (5): we have
$$u=v\circ F^{-1}=v\circ\psi\circ\varphi^{-1}=x_3\circ\varphi^{-1}.$$
 Since $x_3$ is harmonic and
$\varphi$ is biholomorphic, $u$ is a harmonic function.
Differentiating $u(\varphi(z))=x_3(z)$, we obtain
$$2u_z(\varphi(z))\times (-g(z)\omega)=2\frac{\partial x_3}{\partial z}dz=\omega.$$
Hence
$$2u_z(\varphi(z))=\frac{-1}{g(z)}$$
which implies that $||\nabla u||<1$ in $\Omega$ and $||\nabla u||=1$
on $\partial\Omega$.
\cqfd
\section{Proof of Theorem \ref{theorem-periodic}}
\label{appendix3}
We need to construct a meromorphic function $g$ and a holomorphic differential
$\omega$ on a domain $\Sigma\subset\C$ such that $|g|>1$ in $\Sigma$, $|g|=1$ on
$\partial\Sigma$ and $\omega$ is imaginary along $\partial\Sigma$.
Then we proceed as in Section \ref{minimal-vortex} for the correspondence 
``minimal $\to$ vortex''.
Everything depends on the small parameter $t>0$.
\subsection{The domain $\Sigma_t$ and the function $g_t$}
Consider the function
$$f(z)=
\displaystyle c_0+ \sum_{i=1}^n \frac{a_i z}{z-p_i}
.$$
Here $a_1,\cdots,a_n$ are non-zero complex parameters such that
\begin{equation}
\label{normalisation-ai}
a_1+\cdots+a_n=\left\{\begin{array}{ll}0 &\mbox{ in Case (a)}\\
-2c_0 &\mbox{ in Case (b)}\end{array}\right.
\end{equation}
For $t>0$,
let $\Sigma_t$ be the domain $|tf(z)|<1$.
For $t$ small enough, $\Sigma_t$ has $n$ boundary components which
we label $\gamma_1,\cdots,\gamma_n$.
We define the meromorphic
function $g_t$ on $\Sigma_t$ by
$$g_t(z)=\frac{-\ii}{tf(z)}.$$
We have $|g_t|>1$ in $\Sigma_t$ and $|g_t|=1$ on $\partial\Sigma_t$.
\subsection{Opening nodes}
To define the holomorphic differential $\omega_t$ we need to construct the ``double''
of the domain $\Sigma_t$. We do this by ``opening nodes''.
Consider two copies of the complex plane, denoted $\C_1$ and $\C_2$.
As $f$ has a simple pole at $p_i$, there exists a neighborhood $V_i\subset\C_1$ of
$p_i$, a neighborhood $W_i\subset\C_2$ of $\overline{p_i}$ and $\varepsilon>0$ such that
the holomorphic functions
$$v_i(z)=\frac{1}{f(z)}:V_i\to D(0,\varepsilon) \quad\mbox{ and } \quad
w_i(z)=\frac{1}{\overline{f(\overline{z})}}:W_i\to D(0,\varepsilon)$$
are biholomorphic.
Consider the disjoint union $\C_1\cup\C_2$.
For $i=1,\cdots,n$, remove the disks $|v_i|<\frac{t^2}{\varepsilon}$
and $|w_i|<\frac{t^2}{\varepsilon}$. 
Identify the point $z\in V_i$ with the
point $z'\in W_i$ such that
\begin{equation}
\label{opening-nodes}
v_i(z)w_i(z')=t^2.
\end{equation}
This defines a Riemann surface of genus $n-1$ which we denote $\whSigma_t$.
We also define $\Sigma_0$ as the (singular) Riemann surface with $n$ nodes (or double points) obtained by
identifying $p_i\in\C_1$ with $\overline{p_i}\in\C_2$ for $1\leq i\leq n$.
\medskip

The anti-holomorphic involution $\sigma$ which exchanges $z\in\C_1$ with $\overline{z}\in\C_2$ is well defined on $\whSigma_t$ by the following computation:
$$z\sim z'\Rightarrow v_i(z)w_i(z')=t^2\Rightarrow
\overline{v_i(z)}\overline{w_i(z')}=t^2\Rightarrow
w_i(\overline{z})v_i(\overline{z'})=t^2\Rightarrow \sigma(z)\sim\sigma(z').$$
We see $\Sigma_t$ as a domain in $\C_1$.
The involution $\sigma$ exchanges the disjoint domains $\Sigma_t\subset\C_1$
and  $\overline{\Sigma_t}\subset\C_2$ and its fixed set is
$\partial\Sigma_t$, indeed:
$$z=\sigma(z)\Leftrightarrow z\sim \overline{z}\Leftrightarrow
v_i(z)w_i(\overline{z})=t^2\Leftrightarrow |tf(z)|^2=1.$$
Hence we can see $\whSigma_t$ as the disjoint union of $\Sigma_t$ and its
mirror image $\overline{\Sigma_t}$ glued along their boundaries by the map $z\mapsto\overline{z}$. In other words, $\whSigma_t$ is the double of $\Sigma_t$. The point of
constructing $\whSigma_t$ by opening nodes is that it will allow us to understand the limit
$t\to 0$.
Finally, we can extend the definition of the function $g_t$ to $\whSigma_t$ by
$$g_t(z)=\left\{\begin{array}{ll}
\displaystyle\frac{-\ii}{tf(z)} =-\ii \,\frac{v_i(z)}{t}& \mbox{ in $\C_1$}\\
\displaystyle-\ii \,t\overline{f(\overline{z})}
=-\ii\,\frac{t}{w_i(z)} & \mbox{ in $\C_2$}
\end{array}\right.$$
The identification \eqref{opening-nodes} implies that $g_t$ is well defined on $\whSigma_t$.
Moreover, $g_t$ has the symmetry $g_t\circ\sigma=1/\overline{g}$.
\subsection{The holomorphic differential $\omega_t$}
We compactify $\whSigma_t$ by adding the points at infinity in
$\C_1$ and $\C_2$, denoted $\infty_1$ and $\infty_2$.
We also denote $0_1$ and $0_2$ the points $z=0$ in $\C_1$ and $\C_2$.
\begin{proposition}
\label{proposition-omega}
For $t\neq 0$, there exists a unique meromorphic differential $\omega_t$ on $\whSigma_t$ with 4 simple poles at $0_1$, $0_2$, $\infty_1$
and $\infty_2$, such that
\begin{eqnarray*}
&&\int_{\gamma_j}\omega_t=-2\pi \ii c_j\quad\mbox{ for } 1\leq j\leq n\\
&&\Res_{0_1}\omega_t=c_0\\
&&\Res_{0_2}\omega_t=-\overline{c_0}.
\end{eqnarray*}
Here $c_0,c_1,\cdots,c_n$ are given from the configuration. Moreover:
\begin{enumerate}
\item $\omega_t$ has the following symmetry: $\sigma^*\omega_t=-\overline{\omega_t}$.
\item $\omega_t$ extends analytically
at $t=0$ with
\begin{equation}
\label{omega0}
\omega_0=\frac{ c_0\,dz}{z}+\sum_{i=1}^n \frac{c_i\,dz}{z-p_i}  \quad\mbox{ in $\C_1$}
\end{equation}
\item The residues of $g_t\omega_t$ at $0_1$ and $\infty_1$ are given by
\begin{equation}
\label{residu0}
\Res_{0_1}g_t\omega_t=\frac{-\ii}{t}
\end{equation}
\begin{equation}
\label{residu-infty}
\Res_{\infty_1}g_t\omega_t=\frac{\ii}{t}
\end{equation}
\end{enumerate}
\end{proposition}
Proof:
the existence of $\omega_t$ follows from the standard theory of compact Riemann
surface: the curves $\gamma_1,\cdots,\gamma_{n-1}$ are the $A$-cycles of
a canonical homology basis. In general, one can define a meromorphic differential
with simple poles by prescribing its periods on these cycles and the residues at the poles, with
the only restriction that the sum of the residues is zero. In our case, prescribing
the residue at $\infty_1$ is the same as prescribing the period along the last cycle
$\gamma_n$.
\medskip

Proof of Point (1):
$$\int_{\gamma_j}\overline{\sigma^*\omega_t}
=\int_{\sigma(\gamma_j)}
\overline{\omega_t}
=-\int_{\gamma_j}\overline{\omega_t}
=2\pi\ii c_j,$$
$$\Res_{0_2}\overline{\sigma^*\omega_t}=\overline{\Res_{0_1}\omega_t}=\overline{c_0},
\qquad
\Res_{0_1}\overline{\sigma^*\omega_t}=-c_0.$$
Hence the meromorphic differentials $\overline{\sigma^*\omega_t}$ and
$-\omega_t$ have the same poles, periods and residues, so they are equal.
\medskip

Proof of Point (2): we know from the theory of opening nodes that $\omega_t$ extends analytically at $t=0$, and $\omega_0$ is a meromorphic differential on $\Sigma_0$ with
at most simples poles at the nodes
(see \cite{masur}, \cite{nosym} or \cite{crelle}). The residues at the poles are determined
by the prescribed periods. (Here $\gamma_i$ is oriented as a boundary of $\Sigma_t$ so has clockwise orientation.)
\medskip

Proof of Point (3): we have
$$f(0)=c_0,
\quad g_t(0_1)=\frac{-\ii}{tc_0},
\quad\Res_0\omega_t=c_0.$$
This gives \eqref{residu0}.
In Case (a), we have, using $a_1+\cdots+a_n=c_1+\cdots+c_n=0$,
$$f(\infty)=c_0,
\quad g_t(\infty_1)=\frac{-\ii}{tc_0},
\quad\Res_{\infty_1}\omega_t=-c_0.$$
In Case (b), we have, using $a_1+\cdots+a_n=c_1+\cdots+c_n=-2c_0$,
$$f(\infty)=-c_0,
\quad g_t(\infty_1)=\frac{\ii}{tc_0},
\quad \Res_{\infty_1}\omega_t=c_0.$$
In either cases, this gives \eqref{residu-infty}.\cqfd
\subsection{Zeros of $\omega_t$}
\begin{proposition}
\label{proposition-zeros}
For $t$ small enough, one can adjust the parameters
$a_1,\cdots,a_n$ so that $\omega_t$ has a zero at each zero and pole of $g_t$.
Moreover, $a_i(t)$ is a smooth function of $t$ and $a_i(0)=c_i$.
\end{proposition}
Proof:
the meromorphic differential
$\omega_t$ has 4 poles on a genus $n-1$ compact Riemann surface so has $2n$ zeros. By symmetry, $\omega_t$ has $n$ zeros in $\C_1$, which we call
$\zeta_1(t),\cdots,\zeta_n(t)$. We have to solve $f(\zeta_i(t))=0$ for $1\leq i\leq n$.
We use the implicit function theorem at $t=0$.
When $t=0$, an obvious solution is to take $a_i=c_i$ which gives by \eqref{omega0}
\begin{equation}
\label{omega00}
\omega_0=f(z)\frac{dz}{z}.
\end{equation}
Then we compute
$$\frac{\partial f(\zeta_i)}{\partial a_j}|_{t=0}=\frac{\zeta_i}{\zeta_i-p_j}.$$
The determinant of this $n\times n$ matrix is a Cauchy determinant so 
it is invertible.
The problem to apply the implicit function theorem is
that the parameters $a_1,\cdots,a_n$ are constrained by Equation
\eqref{normalisation-ai}, so we have in fact only $n-1$ parameters available.
Let $h(a_1,\cdots,a_n,t)=(f(\zeta_1),\cdots,f(\zeta_{n-1}))$.
\begin{lemma}
\label{lemma1}
The partial differential of $h$ with respect to $(a_1,\cdots,a_n)$ at
$(c_1,\cdots,c_n,0)$, restricted to the space $a_1+\cdots+a_n=0$, is
an isomorphism.
\end{lemma}
Proof: let $(a_1,\cdots,a_n)$ be in the kernel of the partial differential of $h$. Then
$$\sum_{j=1}^n\frac{a_i}{\zeta_i-p_j}=0\quad\mbox{ for }1\leq i\leq n-1.$$
Consider the meromorphic differential on $\C\cup\{\infty\}$
$$\mu=\sum_{j=1}^n\frac{a_i}{z-p_j}\,dz.$$
Then $\mu$ has $n-1$ zeros at $\zeta_1,\cdots,\zeta_{n-1}$, $n$ poles
at $p_1,\cdots,p_n$ and is holomorphic
at $\infty$ because $a_1+\cdots+a_n=0$. Hence $\mu=0$ so $a_1=\cdots=a_n=0$.\cqfd
\medskip

By Lemma \ref{lemma1} and the implicit function theorem, we can solve $f(\zeta_i)=0$ for $1\leq i\leq n-1$.
It  remains to understand why $f(\zeta_n)=0$.
Let $\zeta'$ be the last zero of $f$ in $\C_1$. The meromorphic differential
$g_t\omega_t$ has 4 simple poles at $0_1$, $\infty_1$, $0_2$, $\infty_2$,
and at most a simple pole at $\zeta'$. By Point (3) of Proposition
\ref{proposition-omega} and by symmetry, the sum of the residues of
$g_t\omega_t$ at $0_1$, $\infty_1$, $0_2$ and $\infty_2$ is zero.
By the residue theorem, the residue at $\zeta'$ is zero, so $g_t\omega_t$
is actually holomorphic at $\zeta'$, which means that $\zeta'=\zeta_n$.
This proves Proposition \ref{proposition-zeros}.
\cqfd

\begin{remark} We tacitly assumed that the zeros $\zeta_1,\cdots,\zeta_n$
are distinct, which is the generic case. In case there are multiple zeros, the
proof must be fixed using the Weierstrass preparation theorem. See details
in \cite{nosym}.
\end{remark}
\subsection{The period problem}
From now on, we assume that $a_i(t)$ has the value given by Proposition 
\ref{proposition-zeros}.
\begin{proposition}
\label{proposition-periods}
For $t$ small enough, one can adjust $p_1,\cdots,p_n$ so that
for $1\leq i\leq n$,
\begin{equation}
\label{period-problem}
\int_{\gamma_i}g_t\omega_t=0.
\end{equation}
Moreover, $p_i(t)$ is a smooth function of $t$ and $p_i(0)$ is given by the configuration.
\end{proposition}
Proof: by Point (3) of Proposition \ref{proposition-omega} and the residue theorem , we have
\begin{equation}
\label{sum-periods}
\sum_{i=1}^n \int_{\gamma_i}g_t\omega_t=2\pi\ii \left(\Res_{0_1}g_t\omega_t+
\Res_{\infty_1}g_t\omega_t\right)=0.
\end{equation}
So it suffices to solve \eqref{period-problem} for $1\leq i\leq n-1$.
By symmetry and definition of $g_t$,
$$\int_{\gamma_i}g_t\omega_t
=\int_{\gamma_i}\frac{1}{\overline{g_t}}(-\overline{\omega_t})
=\ii t \int_{\gamma_i}\overline{f\omega_t}.$$
So we want to solve
\begin{equation}
\label{period-problem2}
\int_{\gamma_i}f\omega_t=0\quad \mbox{ for } 1\leq i\leq n-1.
\end{equation}
We solve \eqref{period-problem2} using the implicit function theorem at $t=0$.
We compute
\begin{eqnarray*}
\int_{\gamma_i}f\omega_0
&=&\int_{\gamma_i} f^2(z)\frac{dz}{z}\qquad \mbox{ using \eqref{omega00}}\\
&=&2\pi\ii \Res_{p_i}\left(c_0+\sum_{j=1}^n\frac{c_j z}{z-p_j}\right)^2\frac{dz}{z}\\
&=&2\pi\ii \Res_{p_i}\left[\frac{c_i^2z}{(z-p_i)^2}+\frac{2c_i}{z-p_i}\left(c_0+\sum_{j\neq i}
\frac{c_j z}{z-p_j}\right)\right]\\
&=&2\pi\ii \left(c_i^2+2c_0c_i+2\sum_{j\neq i}\frac{c_ic_jp_i}{p_i-p_j}\right)\\
&=&2\pi\ii\left[c_i^2+2c_0c_i+\sum_{j\neq i}c_ic_j\left(\frac{p_i+p_j}{p_i-p_j}+1\right)\right]\\
&=&2\pi\ii \left(\sum_{j=1}^n c_ic_j+2c_0c_i+\sum_{j\neq i}c_ic_j\frac{p_i+p_j}{p_i-p_j}\right)\\
&=& 2\pi\ii F_i
\end{eqnarray*}
where $F_i$ is as in Definition \ref{definition-periodic-configuration}.
Since the configuration is balanced and non-degenerate, we can solve
\eqref{period-problem2} using the implicit function theorem.
\cqfd
\begin{remark} Since $F_1+\cdots+F_n=0$, it was crucial to have the relation
\eqref{sum-periods} amongst the periods for all $t$.
\end{remark}
\subsection{Solution of the hollow vortex problem}
We proceed as in Appendix \ref{appendix2}, with $\omega$ replaced by $t\,\omega_t$.
We define $\varphi_t$ on $\Sigma_t\subset\C_1$ for $t\neq 0$ by
$$\varphi_t(z)=\ii\log z_0-\int_{z_0}^z g_t \,t\,\omega_t.$$
\begin{proposition}
\label{proposition-final}
\begin{enumerate}
\item $\varphi_t:\Sigma_t\to\C/2\pi\Z$ is a well defined holomorphic map.
\item $\varphi_t$ extends analytically at $t=0$ (away from the points $p_1,\cdots,p_n$) with $\varphi_0(z)=\ii\log z$.
\item $\varphi_t$ is a diffeomorphism onto its image $\Omega_t=\varphi_t(\Sigma_t)\subset\C/2\pi\Z$. (The domain $\Omega_t$ lifts to a periodic domain in the plane with period $2\pi$.)
\item The function $u_t$ defined on $\varphi_t(\Sigma_t)$ by
$$u_t(\varphi_t(z))=\Re\int_{z_0}^zt\,\omega_t$$
solves Problem \eqref{HVP} on $\Omega_t$.
Moreover, the velocity $\vec{v}_t$ and its circulation are given by
\begin{equation}
\label{velocity}
\vec{v}_t(\varphi_t(z))=t\overline{f(z)}
\end{equation}
$$\boC(\varphi_t(\gamma_i))=2\pi\, t\, c_i.$$
\end{enumerate}
\end{proposition}
Proof:
by Propositions \ref{proposition-zeros}, $g_t\omega_t$ is holomorphic and
non-zero in $\Sigma_t$. By Proposition \ref{proposition-periods}, the
only periods of $g_t\omega_t$ come from the residues at $0$ and $\infty$.
By Point (3) of Proposition \ref{proposition-omega}, $\varphi_t$ is well defined modulo $2\pi$, which proves
Point (1). To prove Point (2), we write
$$\varphi_t(z)=\ii\log z_0+\ii\int_{z_0}^z \frac{\omega_t}{f}$$
and we use  Equation \eqref{omega00}.
Regarding Point (3), we already know that $\varphi_t$ is a local diffeomorphism
because its derivative does not vanish. Consider a component $\gamma_i$
of $\partial\Sigma_t$. The unit normal to $\varphi_t(\gamma_i)$ is $g_t$.
Since the function $g_t$ is a diffeomorphism from $\gamma_i$ to the unit circle,
$\varphi_t(\gamma_i)$ is a small convex curve. Consider then the well-defined 
holomorphic function
$\psi=\exp(\ii\varphi_t):\Sigma_t\to\C$. Observe that $\psi$ has a simple pole at
$0$ and a simple zero at $\infty$. Since each $\psi(\gamma_i)$ bounds a disk in the
Riemann sphere, we can extend $\psi$ into a local homeomorphism $\widetilde{\psi}$ from
the Riemann sphere to itself. Since the Riemann sphere is compact and simply connected,
$\widetilde{\psi}$ must be a diffeomorphism. Hence $\varphi_t$ is injective,
which proves Point (3).

\medskip

Proof of Point (4): As in the proof of Point (5) of Proposition \ref{correspondence2},
we have
$$2\frac{\partial u_t}{\partial z}(\varphi_t(z))=\frac{-1}{g_t(z)}=-\ii \,t f(z).$$
Formula \eqref{velocity} for the velocity follows.
From the definition of $u_t$ we obtain
$$\varphi_t^*(2u_z\,dz)=t\omega_t.$$
By \eqref{circulation}, the residue theorem and the definition
of $\omega_t$, we have
$$\boC(\varphi_t(\gamma_i))=\ii\int_{\varphi_t(\gamma_i)}2u_z\,dz
=\ii\int_{\gamma_i}\varphi_t^*(2u_z\,dz)
=\ii\int_{\gamma_i}t\omega_t
=2\pi\,t\, c_i.$$
\cqfd

\bigskip
\noindent
{\sc Martin Traizet\\
Laboratoire de Math\'ematiques et Physique Th\'eorique\\
Universit\'e Fran\c cois Rabelais\\
37200 Tours, France.\\}
{\em email address: }\verb$martin.traizet@lmpt.univ-tours.fr$


\begin{thebibliography}{9}
\bibitem{baker} G. R. Baker, P. G. Saffman, J. S. Sheffield: Structure of a linear array
of hollow vortices of finite cross-section. {J. Fluid Mech.} 74, part 3,
469--476 (1976).
\bibitem{connor} P. Connor, M. Weber:
The construction of doubly periodic minimal surfaces via balance equations.
{\em  American Journal of Mathematics} 134 (5), 1275--1301 (2012).
\bibitem{crowdy} D. Crowdy, C. Green: Analytical solutions for von Karman streets
of hollow vortices. {\em Physics of Fluids} 23, 126602 (2011).
\bibitem{crowdy2} D. Crowdy, J.Roenby:
Hollow vortices, capillary water waves and double quadrature domains.
{\em Fluid Dyn. Res.} 46, 031424 (2014).
\bibitem{eremenko} A. Eremenko, E. Lundberg:
Quasi-exceptional domains.
arXiv:1405.7754v1 (2014).
\bibitem{hhp} F. H\'elein, L. Hauswirth, F. Pacard: A note on some overdetermined elliptic problem. {\em Pacific J. Math.} 250, 319--334 (2011).
\bibitem{hoffman-meeks} D. Hoffman, W.H. Meeks III: The strong halfspace theorem for minimal surfaces. {\em Invent. Math.} 101, 373--377 (1990).
\bibitem{fuchs} W.H.J. Fuchs: A Phragmen Lindel\"of Theorem conjectured by D.J. Newman. {\em Trans. Amer. Math. Soc.} 267, No. 1, 285--293 (1981).
\bibitem{karcher} H. Karcher: The Triply Periodic Minimal Surfaces of {A}lan {S}choen and their Constant Mean Curvature Companions.
{\em Manuscripta Math.} 64, 291--357 (1989).
\bibitem{khavinson} D. Khavinson, E. Lundberg, R. Teodorescu: An overdetermined problem in potential theory.
{\em Pacific J. Math.} 265, 85--111 (2013).
\bibitem{masur}H. Masur: The extension of the Weil-Petersson
metric to the boundary of Teichmuller space. {\em Duke Math. J.} 43, 623--635 (1976).
\bibitem{mazet} L. Mazet: Quelques r\'esultats de non-existence pour l'\'equation des surfaces minimales. {\em Bull. Sci. Math.} 128, 577--586  (2004).
\bibitem {nosym} M. Traizet: An embedded minimal surface with no symmetries.
{\em Journal of Diff. Geom.} 60, 103--153 (2002).
\bibitem{GAFA} M. Traizet: Classification of the solutions to an overdetermined elliptic problem in the plane. {\em Geometric and Functional Analysis} 24 (2), 690--720 (2014).
\bibitem{crelle} Martin Traizet: Opening infinitely many nodes. {\em J. reine angew. Math.} 684, 165--186 (2013).
\bibitem{vandyke} M. Van Dyke: An Album of Fluid Motion. The Parabolic Press, Stanford
University (1982).
\bibitem{weber-wolf} M. Weber, M. Wolf. Teichm\"uller theory and handle addition for minimal surfaces. {\em Annals of Math.} 156, 713--795 (2002). 
\end{thebibliography}
\end{document}